\documentclass[11pt]{article}
\usepackage[utf8]{inputenc}
\usepackage{amsmath, amsthm, amssymb}
\usepackage{mathrsfs}
\usepackage{enumitem}
\usepackage{bm}
\usepackage{tikz}
\usepackage{xcolor}
\usepackage{xypic}
\usepackage{verbatim}
\usepackage{cite}
\usepackage{hyperref}

\usetikzlibrary{quotes,angles}

\hypersetup{
    colorlinks=true,
    linkcolor=blue,
    filecolor=magenta,
    urlcolor=cyan,
}

\theoremstyle{definition}
\newtheorem{theorem}{Theorem}[section]

\newtheorem{corollary}[theorem]{Corollary}
\newtheorem{definition}[theorem]{Definition}

\newtheorem{lemma}[theorem]{Lemma}
\newtheorem{proposition}[theorem]{Proposition}
\newtheorem{remark}[theorem]{Remark}
\newtheorem*{assumption*}{Assumption}

\newcommand{\bkt}[1]{\biggl( #1 \biggr) }
\newcommand{\norm}[1]{\bigl|\bigl| #1 \bigr| \bigr|}
\newcommand{\abs}[1]{\biggl| #1 \biggr| }
\newcommand{\sabs}[1]{\bigl| #1 \bigr| }

\title{Matrix Deviation Inequality for $\ell_{p}$-Norm}

\author{
  Yuan-Chung Sheu
  \thanks{Department of Applied mathematics, National Yang Ming Chiao Tung University, Hsinchu, Taiwan.}
  \thanks{Email: \href{mailto:sheu@math.nctu.edu.tw}{sheu@math.nctu.edu.tw}}
  \and
  Te-Chun Wang
  \thanks{Department of Mathematics and Statistics, University of Victoria, British Columbia, Canada.}
  \thanks{Email: \href{mailto:lieb.am07g@nctu.edu.tw}{lieb.am07g@gmail.com}}
}

\begin{document}

\maketitle

\date{}

\begin{abstract}
Motivated by the general matrix deviation inequality for i.i.d ensemble Gaussian matrix \cite[Theorem 11.1.5]{V}, we show that this property holds for $\ell_{p}$-norm with $1\leq p< \infty$ and i.i.d ensemble sub-Gaussian random matrices, which is a random matrix with i.i.d mean-zero, unit variance, sub-Gaussian entries. As a consequence of our result, we establish the Johnson–Lindenstrauss lemma from $\ell_{2}^{n}$-space to $\ell_{p}^{m}$-space for all i.i.d ensemble sub-Gaussian random matrices.
\end{abstract}

\tableofcontents

\section{Introduction} \label{Introduction}
Given an $m \times n$ random matrix $A$, the uniform deviation inequality plays an important role in theory of random matrices. Also it has many interesting and important consequences. We first quote a classical result \cite[Theorem 11.1.5]{V} for i.i.d ensemble Gaussian random matrices with respect to positive-homogeneous and
subadditive function, which is very useful in asymptotic geometric analysis, and the proof goes back to \cite[Lemma 3]{Schechtman}, which has a different formulation.
\begin{theorem} \label{MDI0}
Let $A \in \mathbb{R}^{m \times  n}$ be a random matrix with i.i.d $\mathscr{N}(0,1)$ entries, $T \subseteq \mathbb{R}^{n}$, and $f:\mathbb{R}^{m} \mapsto \mathbb{R}$ be a function such that 
\begin{equation*}
    f(c v) = c f(v) \quad \forall c \geq 0, \; v \in \mathbb{R}^{m} \quad \text{and}\quad
    f(u+v) \leq f(u) + f(v) \quad \forall u,v \in \mathbb{R}^{m}.
\end{equation*}
Then we have
\begin{equation} \label{MDI}
    \mathbb{E}[\sup_{x \in T}\sabs{f(Ax) - \mathbb{E}[f(Ax)]}] \leq C \emph{Lip}(f) \gamma(T),
\end{equation}
where $\gamma(T) \equiv \mathbb{E}[\sup_{x\in T} |\langle g,x \rangle|]$, $g \sim \mathscr{N}(0,I_{n})$, and $\emph{Lip}(f)$ is the Lipschitz constant of $f$. Here $C$ is an absolute universal constant and $\gamma(T)$ is called the Gaussian complexity of $T$.
\end{theorem}
Note that if $f(x) \equiv \sup_{y\in S}\langle x, y \rangle$, then $f(x)$ satisfies the conditions of Theorem \ref{MDI0} and
$\text{Lip}(f) = \text{rad}(S)$, so (\ref{MDI}) is sharp (see \cite[Theorem 11.2.4]{V} and \cite[Exercise 8.7.2]{V}). An important application of matrix deviation inequality is to establish a Lipschitz embedding between two normed spaces. See \cite[Theorem 7]{Schechtman} and \cite[Theorem 11.3.3]{V} for using i.i.d ensemble Gaussian matrix to prove the existence of the embedding from a finite dimensional normed spaces into a low-dimensional Euclidean spaces and the embedding from a high-dimensional Euclidean spaces into a low-dimensional
Euclidean spaces. See also \cite[Theorem 8.7.1]{V} and \cite{Tstable} for Chevet inequality and \cite{Chen} for the $\ell_{p}$-Gaussian-Grothendieck problem. 

Still, it is a challenging open problem to study the universality of the matrix deviation inequality (see \cite[Remark 11.1.9]{V}). In other words, whether the general matrix deviation inequality holds for i.i.d ensemble sub-Gaussian random matrix, which is a random matrix with i.i.d mean-zero, unit variance, sub-Gaussian entries. Note that this problem was solved by Liaw et al. \cite[Theorem 1]{V2} when $f(\cdot) = ||\cdot||_{2}$. Namely, if $A \in \mathbb{R}^{m \times n}$ is an i.i.d ensemble sub-Gaussian random matrix with $K = ||A_{1,1}||_{\psi_{2}}$, where $||X||_{\psi_{2}}$ is the sub-Gaussian norm of $X$ (see Definition \ref{sub-gaussian}), then
\begin{equation} \label{origin}
\norm{\sup_{x \in T}\sabs{||Ax||_{2} - m^{1/2}||x||_{2}}\;}_{\psi_{2}}
\leq
C K^{2}\gamma(T) \quad \forall T \subseteq \mathbb{R}^{n}.
\end{equation}
In addition, (\ref{origin}) still holds if $m^{1/2}||x||_{2}$ is replaced by $\mathbb{E}[||Ax||_{2}]$. See \cite[Chapter 9]{V} for the applications. 

\paragraph{Main Results.}
In this paper, we aim to prove the matrix deviation inequality for $\ell_{p}$-norm, $1\leq p < \infty$, and i.i.d ensemble sub-Gaussian random matrices.
\begin{theorem} \label{main}
Let $A \in \mathbb{R}^{m \times n}$ be a random matrix with i.i.d, mean-zero, unit variance, sub-Gaussian entries $\{A_{i,j}\}$ and $K = ||A_{i,j}||_{\psi_{2}}$. Then we have
\begin{align} \label{goal2}
    \norm{\sup_{x \in T}\sabs{||Ax||_{p} - & m^{
    \frac{1}{p}}||A_{1}x||_{L^{p}}}\;}_{\psi_{2}}\nonumber \\
    &\leq
    \begin{cases}
    C_{p} K^{p} ||A_{1,1}||_{L^{p}}^{-(p+2)} \emph{Lip}(||\cdot||_{p})\gamma(T), &\text{if } p\in [1,2]\\
    C_{p} K^{p+2}   \emph{Lip}(||\cdot||_{p})\gamma(T), &\text{if } p\in(2,\infty).
    \end{cases}
\end{align}
where $A_{i}$ is the ith row of $A$, $||A_{1}x||_{L^{p}} = \mathbb{E}[|A_{1}x|^{p}]^{\frac{1}{p}}$, and $C_{p}$ is a positive absolute constant depending only on $p$. In addition, (\ref{goal2}) still holds if $m^{1/p}||A_{1}x||_{L^{p}}$ is replaced by $\mathbb{E}[||Ax||_{p}]$.
\end{theorem} 
As a consequence of Theorem \ref{main}, we show that any i.i.d ensemble sub-Gaussian random matrix can be regraded as an embedding from $\ell_{2}^{n}$-space to $\ell^{m}_{p}$-space such that the distances between two points do not increase by more than a factor $D_{p}(1+\epsilon)$ and do not decrease by more than a factor $d_{p}(1-\epsilon)$. 
\begin{corollary} \label{JL}
Let $\epsilon \in (0,1)$, $T$ be a finite subset of $\mathbb{R}^{n}$ containing $N$ elements, and
\begin{equation} \label{d_and_D}
    (d_{p},D_{p}) =
    \begin{cases}
    (C_{p}||A_{1,1}||_{p},1), & \text{ if } 1\leq p < 2,\\
    (1,1), &\text{ if } p = 2,\\
    (1,C_{p}'K), &\text{ if } 2 < p < \infty.
    \end{cases}
\end{equation}
Then, under assumption of Theorem \ref{main}, we have
\begin{align} \label{JL_EQ}
    &\mathbb{P}\bkt{
    d_{p}(1-\epsilon)||x-y||_{2} \leq \norm{\frac{1}{m^{\frac{1}{p}}}A(x-y)}_{p} \leq D_{p}(1+\epsilon)||x-y||_{2} \quad \forall x,y \in T} \nonumber\\
    &\geq 
    \begin{cases}
    1-2\exp(-C_{p}\frac{\epsilon^{2} m}{K^{2p} ||A_{1,1}||_{p}^{-2(p+3)}\log(N)}), &\text{ if } 1\leq p\leq 2\\
    1-2\exp(-C_{p}\frac{\epsilon^{2} m^{2/p}}{K^{2(p+2)}\log(N)}), &\text{ if } 2<p<\infty.
    \end{cases}
\end{align}
\end{corollary}
\begin{remark}
For the problem of dimension reduction, Brinkman and Charikar \cite{intro1} and Ping Li \cite{intro2} both give overviews of the results in this area. See also \cite{diamond} for the problem of distortion and \cite[Section 11.3]{V} for the random projection. For more general results and similar problems, see \cite{sJ3}, \cite{sJ2}, and \cite{sJ1}.
\end{remark}

\paragraph{Heuristics.} The core of our proof of Theorem \ref{main} is to show that 
\begin{equation*}
    R_{x} \equiv ||Ax||_{p} - m^{\frac{1}{p}}||A_{1}x||_{L^{p}} 
    \quad  \forall x\in \mathbb{R}^{n}
\end{equation*}
has sub-Gaussian increments (see Lemma \ref{v2case3}). To do this, we will use the approach given by \cite[The proof of Lemma 3]{Schechtman} and \cite[p.292]{V2}, which indicates that it suffices to consider some special cases of Lemma \ref{v2case3} (see Lemma \ref{v2case1} and Lemma \ref{v2case2}). In order to prove these special cases, we establish the sub-Gaussian concentration inequality with respect to $\ell_{p}$-norm by using \cite[
Corollary 1.4]{a_sub_exponential} to control the tail probability of the sum of i.i.d $2/p$-Orlicz random variables $|X_{j}|^{p}$ (see Lemma \ref{v2main}), where $X_{1}$ is a mean-zero sub-Gaussian random variable (note that Lemma \ref{v2main} is a generalization of \cite[Proposition 5.1]{V2}).

\section{$\alpha$-Orlicz Random Variables} \label{Orlicz}
\begin{definition} \label{sub-gaussian}
Let $\alpha > 0$ and $X$ be a random variable. The $\alpha$-Orlicz norm of $X$ is defined by
\begin{equation} \label{definition}
    ||X||_{\psi_{\alpha}} \equiv \inf\{t>0 : \mathbb{E}[\exp(|X|^{\alpha}/t^{\alpha})]\leq 2\}
\end{equation}
(For convenience, we set $\inf \emptyset = \infty$).
\end{definition}
We say $X$ is a $\alpha$-Orlicz random variable if $||X||_{\psi_{\alpha}} < \infty$.
In particular, we say $X$ is a sub-Gaussian random variable if $||X||_{\psi_{2}} < \infty$. Note that $||\cdot||_{\psi_{\alpha}}$ is a norm if and only if $\alpha \geq 1$. Nevertheless, $||\cdot||_{\psi_{\alpha}}$ still make sense for any $\alpha > 0$. 

The following proposition states the equivalent definitions of $||\cdot||_{\psi_{\alpha}}$, which will be used throughout this paper. Note that the proof of Proposition \ref{intro1} is the same as \cite[Proposition 2.5.2]{V}. 

\begin{proposition} \label{intro1}
Let $\alpha > 0$ and $X$ be a random variable. Then the following properties are equivalent:
\begin{enumerate} [label=(\alph*)]
  \item The MGF of $|X|^{\alpha}$ is bounded at some point, namely
  \begin{equation} \label{1p}
  \mathbb{E}[\exp(|X|^{\alpha}/K_{1}^{\alpha})] \leq 2.
  \end{equation}
  \item The tails of $X$ satisfy
  \begin{equation} \label{2p}
  \mathbb{P}(|X| \geq t) \leq 2\exp(-t^{\alpha}/K_{2}^{\alpha}) \quad \forall t\geq 0.
  \end{equation}
  \item The moments of $X$ satisfy
  \begin{equation} \label{3p}
  ||X||_{L^{p}} \leq K_{3}   p^{\frac{1}{\alpha}} \quad \forall p\geq \alpha
  \end{equation}
\end{enumerate}
Here the parameters $K_{i} > 0$ appearing in these properties differ from each other by at most a constant that depends on $\alpha$.
\end{proposition}

By the definition of $||\cdot||_{\psi_{\alpha}}$, it is clear that we have the following relation. 
\begin{lemma} \label{cs}
Let $X$ be a random variable such that $||X||_{\psi_{\alpha \beta}} \vee || \;|X|^{\beta}||_{\psi_{\alpha}} < \infty$. Then $||\;|X|^{\beta}||_{\psi_{\alpha}} =  ||X||_{\psi_{\alpha\beta}}^{\beta}$.
\end{lemma}

\section{Proofs} \label{Main Proof}
\subsection{Proof of Theorem \ref{main}} \label{Main}
Consider the norm induced by $A$ as follows and recall the definition of $R_{x}$:
\begin{equation} \label{randomfield}
||x|| \equiv ||A_{1}x||_{L^{p}} 
\quad 
\text{and}
\quad
R_{x} 
= ||Ax||_{p} - m^{\frac{1}{p}}||x|| \quad  \forall x\in \mathbb{R}^{n}.
\end{equation}
Recall the Generic chaining bound \cite[Theorem 8.5.3]{V}. Note that the proof of \cite[Theorem 8.5.3]{V} actually gives the following estimation.
\begin{proposition} \label{Talagrand}
Let $\mathscr{T} \subseteq \mathbb{R}^{n}$, $x_{0} \in \mathscr{T}$, and $\{\mathscr{R}_{x}\}_{x\in \mathscr{T}}$ be a random process such that $\norm{\mathscr{R}_{x} - \mathscr{R}_{y}}_{\psi_{2}} \leq \mathscr{K} ||x-y||_{2}$
for every $x,y \in \mathscr{T}$. Then
\begin{equation*}
    \norm{\sup_{x\in \mathscr{T}}|\mathscr{R}_{x} - \mathscr{R}_{x_{0}}|}_{\psi_{2}} \leq C \mathscr{K}  \gamma(\mathscr{T}),
\end{equation*}
where $C$ is a positive absolute constant. In particular, if $\mathscr{R}_{x_{0}} = 0$, then
$
\norm{\sup_{x \in \mathscr{T}} |\mathscr{R}_{x}|}_{\psi_{2}} \leq C\mathscr{K}\gamma(\mathscr{T}).
$    
\end{proposition}

To prove Theorem \ref{main}, it suffices to show that $\{R_{x}\}_{x \in \mathbb{R}^{n}}$ has sub-Gaussian increments (i.e., Lemma \ref{v2case3}). Indeed, since $\gamma(T\bigcup \{0\}) = \gamma(T)$ and $R_{0} = 0$, it follows that (\ref{goal2}) is an immediate consequence of Proposition \ref{Talagrand} and Lemma \ref{v2case3}. Also, by triangle inequality and $\sup_{x\in T}||x||_{2} C= \sup_{x\in T} \mathbb{E}[|\langle g,x \rangle|] =\gamma(T)$, where $C = \mathbb{E}[|g|]$ and $g \sim \mathscr{N}(0,1)$, it is clear that $m^{1/p}||A_{1}x||_{L^{p}}$ can be replaced by $\mathbb{E}[||Ax||_{p}]$.

Let us start with some properties that will be used throughout the proof. 
\begin{enumerate} [label=(\alph*)]
    \item
    Applying Jensen's inequality shows that
    \begin{equation} \label{K}
    K
    \geq \inf\{t> 0:\exp(\frac{\mathbb{E}[|A_{1,1}|^{2}]}{t^{2}}) \leq 2\}
    =\sqrt{\frac{1}{\ln2}} > 1.
    \end{equation}
    \item  Note that
    $||\cdot||$ and $||\cdot||_{2}$ are equivalent. Namely, 
    \begin{equation} \label{norm}
    C_{p} ||A_{1,1}||_{L^{p}}||x||_{2} \leq ||x|| \leq ||x||_{2} \quad \forall 1\leq p \leq 2
    \end{equation}
    and
    \begin{equation} \label{norm2}
    ||x||_{2} \leq ||x|| \leq C_{p}'K ||x||_{2} \quad \forall 2\leq p <\infty,
    \end{equation}
    where $C_{p}$ and $C'_{p}$ are positive constants that depend on $p$. The proof of (\ref{norm2}) follows from \cite[Exercise 2.6.5]{V}. The lower bound of (\ref{norm}) is an immediate consequence of Marcinkiewicz–Zygmund inequality \cite[Section 10.3]{chow1997probability} and Minkowski's integral inequality \cite[Theorem 6.2.7]{stroock}. Indeed, if $||x||_{2} = 1$, then
    \begin{equation*}
        ||x|| \geq  C_{p} \mathbb{E}[(\sum_{j=1}^{n} (A_{1,j}x_{j})^{2})^{p/2}]^{1/p} 
        \geq C_{p}(\sum_{j=1}^{n}x_{j}^{2} \mathbb{E}[|A_{1,j}|^{p}]^{1/p})^{1/2}
        = C_{p}||A_{1.1}||_{L^{p}}.
    \end{equation*}
    \item Applying Hölder inequality gives
    \begin{equation} \label{Lip}
        \text{Lip}(||\cdot||_{p}) = 
        \begin{cases}
            m^{\frac{1}{p} - \frac{1}{2}}, &\text{ if } 1\leq p\leq2\\
            1, &\text{ if } 2\leq p<\infty.
        \end{cases}
    \end{equation}
\end{enumerate}

\subsubsection{Case 1: $x \in \mathbb{R}^{n}$ and $y = 0$}
\begin{lemma} \label{v2case1}
Under assumption of Theorem \ref{main}, we have
\begin{equation} \label{v2case1_ineq}
\norm{R_{x}}_{\psi_{2}} 
\leq
\begin{cases}
    C_{p} (K/||A_{1,1}||_{L^{p}})^{p}   \emph{Lip}(||\cdot||_{p})||x||_{2}, &\text{ if } 1\leq p\leq 2\\
    C_{p} K^{p}   \emph{Lip}(||\cdot||_{p})||x||_{2}, &\text{ if } 2< p<\infty.
\end{cases}
\end{equation}
\end{lemma}
To prove Lemma \ref{v2case1}, it suffices to establish Lemma \ref{v2main}. Indeed, by \cite[Proposition 2.6.1]{V}, we have $||A_{1}\frac{x}{||x||}||_{\psi_{2}} \leq CK\frac{||x||_{2}}{||x||}$,
so applying the following lemma gives
\begin{align*}
    &||R_{x}||_{\psi_{2}} \leq C_{p} ||x|| K^{p} (\frac{||x||_{2}}{||x||})^{p} \leq C_{p} K^{p} ||A_{1,1}||_{L^{p}}^{-p} ||x||_{2}\quad \text{if } 1\leq p\leq 2; \\
    &||R_{x}||_{\psi_{2}} \leq
    C_{p} ||x|| K^{p} (\frac{||x||_{2}}{||x||})^{p} = C_{p} ||x||_{2} K^{p} (\frac{||x||_{2}}{||x||})^{p-1} 
    \leq C_{p} ||x||_{2} K^{p} \quad \text{if } 2<p<\infty.
\end{align*}
\begin{lemma} \label{v2main}
Let $1\leq p< \infty$ and $\{X_{i}\}_{1\leq i <
\infty}$ be i.i.d sub-Gaussian random variables such that $||X_{1}||_{L^{p}} = 1$ and $K \equiv ||X_{1}||_{\psi_{2}} < \infty$. Then, for each $m\geq 1$ and $X^{(m)} = (X_{1},...,X_{m})$, we have
\begin{equation*}
   \norm{||X^{(m)}||_{p} - m^{\frac{1}{p}}}_{\psi_{2}}
   \leq C_{p} K^{p} \emph{Lip}(||\cdot||_{p}),
\end{equation*}
where $C_{p}$ is a positive absolute constant depending only on $p$.
\end{lemma}
\begin{proof} To prove Lemma \ref{v2main},
it suffices to show that
\begin{equation} \label{v2_main_1}
    \mathbb{P}\bkt{\abs{||X^{(m)}||_{p} - m^{\frac{1}{p}}} \geq s}
    \leq
    \begin{cases}
        2\exp(-C_{p} \frac{s^{2}}{K^{2p}m^{\frac{2}{p}-1}}), & \text{ if } 1\leq p<2\\
        2\exp(-C_{p} \frac{s^{2}}{K^{2p}}), & \text{ if } 2\leq p<\infty\\
    \end{cases}
    \quad \forall  s>0,
\end{equation}
where $C_{p}$ is a positive absolute constant.
\paragraph{Step 1.} In this step, we prove (\ref{v2_main_1}) when $1\leq p < \infty$ and $s \leq K^{p}m^{\frac{1}{p}}$. Note that if $|z - 1| \geq \delta$ and $z\geq 0$, then $|z^{p} - 1| \geq \delta$. Then we have
$$
\mathbb{P}\bkt{\abs{\frac{1}{m^{\frac{1}{p}}}||X^{(m)}||_{p} - 1} \geq \delta}
\leq \mathbb{P}\bkt{\abs{\frac{1}{m}\sum_{i=1}^{m}(|X_{i}|^{p}-1)} \geq \delta}.
$$
Since $||X_{i}||_{\psi_{p}} \leq c_{p}||X_{i}||_{\psi_{2}}$ for each $1\leq p<2$, we have
\begin{equation} \label{p-1_1}
    \norm{|X_{i}|^{p} -1}_{\psi_{1}}
    \leq C_{p}\norm{|X_{i}|^{p}}_{\psi_{1}}
    \leq C_{p}\norm{X_{i}}_{\psi_{p}}^{p} 
    \leq C_{p} K^{p}
    \quad \forall 1\leq p < 2
\end{equation}
and
\begin{equation} \label{p-1_2}
    \norm{|X_{i}|^{p} -1}_{\psi_{\frac{2}{p}}}
    \leq C_{p}\norm{|X_{i}|^{p}}_{\psi_{\frac{2}{p}}}
    \leq C_{p}\norm{X_{i}}_{\psi_{2}}^{p} 
    \leq C_{p} K^{p} 
    \quad \forall 2\leq p < \infty
\end{equation}
by \cite[Lemma A.3]{a_sub_exponential} and Lemma \ref{cs}. Let $a_{i} = \frac{1}{m}$. Then applying \cite[
Corollary 1.4]{a_sub_exponential} with $\alpha = 1$ if $1\leq p < 2$; $\alpha = \frac{2}{p}$ if $2\leq p<\infty$ gives
\begin{align} \label{just_like_0}
      \mathbb{P}\bkt{\abs{\frac{1}{m}\sum_{i=1}^{m}(|X_{i}|^{p}-1)} \geq \delta}
      &\leq 2\exp(-C_{p}\min\{\frac{\delta^{2}}{K^{2p}},\frac{\delta}{K^{p}}\}m) \nonumber\\
      &= 2\exp(-C_{p}\frac{\delta^{2}m}{K^{2p}})  \quad \forall \delta \leq K^{p}, \quad 1\leq p<2
\end{align}
and
\begin{align} \label{just_like}
      &\mathbb{P}\bkt{\abs{\frac{1}{m}\sum_{i=1}^{m}(|X_{i}|^{p}-1)} \geq \delta}
      \leq 2\exp(-C_{p}\min\{\frac{\delta^{2}m}{K^{2p}},\frac{\delta^{\alpha} m^{\alpha}}{K^{\alpha p}}\})\nonumber \\
      &\leq 2\exp(-C_{p}\min\{\frac{\delta^{2}}{K^{2p}},\frac{\delta^{\alpha} }{K^{\alpha p}}\}m^{\alpha})
      \leq 2\exp(-C_{p}\frac{\delta^{2}m^{\alpha}}{K^{2p}})  \quad \forall \delta \leq K^{p}, \quad 2\leq p<\infty.
\end{align}
Therefore, taking $s = \delta m^{\frac{1}{p}}$ proves (\ref{v2_main_1}) when $1\leq p < \infty$ and $s \leq K^{p}m^{\frac{1}{p}}$.

\paragraph{Step 2.} In this step, we prove (\ref{v2_main_1}) when $1\leq p < 2$ and $s > Km^{\frac{1}{p}}$. In fact, we only need to prove (\ref{v2_main_1}) when $1\leq p < 2$ and $s > Km^{\frac{1}{p}} \xi_{p}$, where $\xi_{p}$ is a positive constant that depends on $p$. Indeed, if $\xi_{p} > 1$ and $Km^{1/p} < s < Km^{1/p} \xi_{p}$, then using the result proved in \textbf{Step 1} gives
\begin{align*}
    &\mathbb{P}\bkt{\abs{||X^{(m)}||_{p} - m^{\frac{1}{p}}} \geq s}
    \leq \mathbb{P}\bkt{\abs{||X^{(m)}||_{p} - m^{\frac{1}{p}}} \geq Km^{1/p}}\\
    &\leq 2\exp(-C_{p} \frac{(Km^{1/p})^{2}}{K^{2p}m^{\frac{2}{p}-1}})
    \leq 2\exp(-\frac{C_{p}}{\xi_{p}^{2}} \frac{s^{2}}{K^{2p}m^{\frac{2}{p}-1}}).
\end{align*}
Note that $|a^{r} - b^{r}| \leq |a-b|^{r}$ if $0 < r \leq 1$ and $a,b > 0$. Hence, the tail probability can be estimated as follows:
\begin{align} 
    &\mathbb{P}\bkt{\abs{\frac{1}{m^{\frac{1}{p}}}||X^{(m)}||_{p} - 1} \geq \delta}
    \leq \mathbb{P}\bkt{\abs{\frac{1}{m} \sum_{i=1}^{m} (|X_{i}|^{p} - 1)} \geq \delta^{p}} \label{decomposition0}\\
    &\leq \mathbb{P}\bkt{\frac{1}{m} \sum_{i=1}^{m} (|X_{i}|^{p} - 1) \geq \delta^{p}}+\mathbb{P}\bkt{\frac{1}{m} \sum_{i=1}^{m} (1-|X_{i}|^{p}) \geq \delta^{p}}. \label{decomposition}
\end{align}
Both of the above terms can be controlled by the same argument. In the following, we only estimate the first term. Note that $(\ref{p-1_2})$ holds for $1\leq p< 2$ as well. Hence, applying \cite[Proposition 5.2]{alpha'} with random variable $|X_{i}|^{p}-1$ and $\alpha = \frac{2}{p}$ gives
\begin{align*}
    &\mathbb{P}\bkt{\frac{1}{m} \sum_{i=1}^{m} (|X_{i}|^{p} - 1) \geq \delta^{p}}
    \leq \exp(-\lambda t + m C_{\alpha}^{\alpha'} K^{p\alpha'} \lambda^{\alpha'}) \quad \forall \lambda \geq \frac{1}{K^{p} C_{\alpha}},
    \end{align*}
where $C_{\alpha}$ is a positive constant that depends on $\alpha$, $t = m\delta^{p}$, $\alpha'$ is the Hölder conjugates of $\alpha$. Note that
\begin{align*}
    &(\frac{t}{m K^{p\alpha'} C_{\alpha}^{\alpha'} })^{\frac{1}{\alpha'-1}} \geq \frac{1}{K^{p} C_{\alpha}}
    \iff \delta \geq K C_{\alpha}^{1/p}.
\end{align*}
Hence, if $\delta \geq K C_{\alpha}^{1/p}$ and $\lambda \equiv
(\frac{t}{m K^{p\alpha'} C_{\alpha}^{\alpha'}})^{\frac{1}{\alpha'-1}}$, then
\begin{align*}
    &\mathbb{P}\bkt{\frac{1}{m} \sum_{i=1}^{m} (|X_{i}|^{p} - 1) \geq \delta^{p}}\\
    &\leq \exp\bkt{-\bkt{\frac{t}{mK^{p\alpha'} C_{\alpha}^{\alpha'}}}^{\frac{1}{\alpha'-1}}t+m K^{p\alpha'} C_{\alpha}^{\alpha'} \bkt{\frac{t}{mK^{p\alpha'} C_{\alpha}^{\alpha'}}}^{\frac{\alpha'}{\alpha'-1}}}
    = \exp\bkt{-C_{p} \frac{\delta^{2} m}{K^{2}}},
\end{align*}
where $C_{p}$ is a constant that depends on $p$. Therefore, we obtain
\begin{equation*}
    \mathbb{P}\bkt{\abs{\frac{1}{m^{\frac{1}{p}}}||X^{(m)}||_{p} - 1} \geq \delta}
    \leq 2\exp(-C_{p} \frac{\delta^{2}m}{K^{2}})
    \leq 2\exp(-C_{p} \frac{\delta^{2}m}{K^{2p}}) \quad \forall \delta \geq C_{\alpha}^{1/p} K
\end{equation*}
by using (\ref{K}), so we complete the proof of (\ref{v2_main_1}) when $1\leq p < 2$ and $s > K m^{\frac{1}{p}}$.

\paragraph{Step 3.} In this step, we prove (\ref{v2_main_1}) when $2 \leq p < \infty$ and $s > m^{\frac{1}{p}}K^{p}$. 
Decompose the tail probability as (\ref{decomposition0}) and apply \cite[Corollary 1.4]{a_sub_exponential}. Then we have
\begin{align*}
    \mathbb{P}\bkt{\abs{\frac{1}{m^{\frac{1}{p}}}||X^{(m)}||_{p} - 1} \geq \delta}
    &\leq 2\exp(-C_{p}\min\{\frac{\delta^{2p}}{K^{2p}},\frac{\delta^{p\alpha} }{K^{\alpha p}}\}m^{\alpha})\\
    &= 2\exp(-C_{p}\frac{m^{\alpha} \delta^{2}}{K^{2}})
    \quad\forall \delta> K^{p},
\end{align*}
so, taking $s = \delta m^{\frac{1}{p}}$, we complete the proof of this step. 
\end{proof}
\subsubsection{Case 2: $||x|| = ||y|| = 1$}
\begin{lemma} \label{v2case2}
Under assumption of Theorem \ref{main}, we have
\begin{align} \label{v2case1_ineq}
    &\norm{R_{x} - R_{y}}_{\psi_{2}} \nonumber \\
    &\leq
    \begin{cases}
    C_{p} (K/||A_{1,1}||_{L^{p}})^{p}   \emph{Lip}(||\cdot||_{p})||x-y||_{2}, &\text{ if } 1\leq p\leq 2\\
    C_{p} K^{p}   \emph{Lip}(||\cdot||_{p})||x-y||_{2}, &\text{ if } 2< p<\infty 
    \end{cases}
    \quad \forall \;||x|| = ||y|| = 1.
\end{align}
\end{lemma}
\begin{proof} To prove Lemma (\ref{v2case2}), it suffices to show that
\begin{equation}\label{proof_case2}
    \mathbb{P} \bkt{\abs{\frac{||Ax||_{p} - ||Ay||_{p}}{||x-y||}} \geq s} 
    \leq
    \begin{cases}
        4\exp(-C_{p}\frac{s^{2}}{(K/||A_{1,1}||_{p})^{2p}m^{\frac{2}{p}-1}}), &\text{ if } 1\leq p\leq 2\\
        4\exp(-C_{p}\frac{s^{2}}{K^{2p}}), &\text{ if } 2< p<\infty,
    \end{cases}  
\end{equation}
where $C_{p}$ is a positive constant that depends $p$. Indeed, since $\mathbb{E}[|Z|^{N}] = \int_{0}^{\infty} N s^{N-1} \mathbb{P}(|Z| \geq s) ds$, it is clear that (\ref{proof_case2}) implies (\ref{v2case1_ineq}). Note that if $s \geq 2m^{\frac{1}{p}}$, then (\ref{proof_case2}) is an immediate consequence of Lemma \ref{v2case1}. Indeed, if $u = \frac{x-y}{||x-y||}$, then
\begin{align*}
    &\mathbb{P}(\abs{\frac{||Ax||_{p} - ||Ay||_{p}}{||x-y||}}\geq s)
    \leq \mathbb{P}(||Au||_{p}\geq s) \nonumber
    = \mathbb{P}(||Au||_{p} - m^{\frac{1}{p}}\geq s-m^{\frac{1}{p}}) \\
    &\leq \mathbb{P}(\abs{||Au||_{p} - m^{\frac{1}{p}}} \geq \frac{s}{2}),
\end{align*}
so applying Lemma \ref{v2case1} gives (\ref{proof_case2}). Thus, it remains to prove (\ref{proof_case2}) when $s < 2m^{\frac{1}{p}}$. Since $a^{p-1}|a-b|\leq|a^{p} - b^{p}|$ if $1\leq p<\infty$ and $a,b > 0$, it follows that
\begin{align*}
    &\mathbb{P}\bkt{\abs{\frac{||Ax||_{p} - ||Ay||_{p}}{||x-y||}} \geq s}
    \leq \mathbb{P}\bkt{\abs{\frac{||Ax||_{p}^{p} - ||Ay||_{p}^{p}}{||x-y||}} \geq s||Ax||_{p}^{p-1}} \\
    &\leq \mathbb{P}\bkt{\abs{\frac{||Ax||_{p}^{p} - ||Ay||_{p}^{p}}{||x-y||}} \geq s||Ax||_{p}^{p-1},||Ax||_{p} \geq \frac{m^{\frac{1}{p}}}{2}} 
    + \mathbb{P}\bkt{||Ax||_{p} < \frac{m^{\frac{1}{p}}}{2}}\\
    &\leq \mathbb{P}\bkt{\abs{\frac{||Ax||_{p}^{p} - ||Ay||_{p}^{p}}{||x-y||}} \geq \frac{sm^{\frac{1}{q}}}{2^{p-1}}} + \mathbb{P}\bkt{||Ax||_{p} < \frac{m^{\frac{1}{p}}}{2}} \equiv \mathscr{A}_{1} + \mathscr{A}_{2}.
\end{align*}
Since $s < 2m^{\frac{1}{p}}$ and $||x|| = 1$, applying Lemma (\ref{v2case1}) gives
\begin{align} \label{A2_New}
    &\mathscr{A}_{2}
    \leq \mathbb{P}\bkt{\abs{||Ax||_{p} -m^{\frac{1}{p}}||x||} \geq
    \frac{m^{\frac{1}{p}}}{2}}
    \leq \mathbb{P}\bkt{\abs{||Ax||_{p} -m^{\frac{1}{p}}||x||} \geq\frac{s}{4}} \nonumber\\
    &\leq 
    \begin{cases}
        2\exp(-C_{p}\frac{s^{2}}{(K/||A_{1,1}||_{L^{p}})^{2p} m^{\frac{2}{p}-1}}), &\text{ if }1\leq p\leq 2\\
        2\exp(-C_{p}\frac{s^{2}}{K^{2p}}), &\text{ if }2< p < \infty.
    \end{cases}
\end{align}
To estimate $\mathscr{A}_{1}$, we write $\mathscr{A}_{1}$ as
\begin{equation*}
    \mathscr{A}_{1}
    = \mathbb{P}(\sabs{\frac{1}{m}\sum_{i=1}^{m} \frac{|A_{i}x|^{p} - |A_{i}y|^{p}}{||x-y||}} \geq \delta), 
    \quad
    \text{where}
    \quad
    \delta \equiv \frac{s}{2^{p-1}m^{\frac{1}{p}}},
\end{equation*} 
so it suffices to show that
\begin{equation} \label{Final}
    \norm{\frac{|A_{i}x|^{p} - |A_{i}y|^{p}}{||x-y||}}_{\psi_{\frac{2}{p}}}
    \leq 
    \begin{cases}
        C_{p} (K/||A_{1,1}||_{p})^{p}, &\text{ if } 1\leq p\leq 2\\
        C_{p} K^{p}, &\text{ if } 2<p<\infty.
    \end{cases}
\end{equation}
Indeed, since $\delta = \frac{s}{2^{p-1}m^{\frac{1}{p}}} \leq 2^{2-p} \leq 2K^{p}$ for every $1\leq p<\infty$, applying \cite[Corollary 1.4]{a_sub_exponential} similar to (\ref{just_like_0}) and (\ref{just_like}) yields
\begin{equation} \label{A1_New}
    \mathscr{A}_{1} \leq 
    \begin{cases}
          2\exp(-C_{p}\frac{(\delta/2)^{2}m}{(K/||A_{1,1}||_{p})^{2p}} )
          = 2\exp(-C_{p}\frac{s^{2}}{(K/||A_{1,1}||_{p})^{2p} m^{\frac{2}{p}-1}}), &\text{ if } 1\leq p \leq 2\\
          2\exp(-C_{p}\frac{(\delta/2)^{2}m^{2/p}}{K^{2p}})
          = 2\exp(-C_{p}\frac{s^{2}}{K^{2p}}), &\text{ if } 2<p<\infty.
    \end{cases}
\end{equation}
Hence, it remains to prove (\ref{Final}). Note that $|a^{p} - b^{p}|\leq p|a-b| \sqrt{a^{2p-2}+b^{2p-2}}$ if $1\leq p<\infty$ and $a,b > 0$. Thus, we have
\begin{equation*}
    \norm{|A_{i}x|^{p} - |A_{i}y|^{p}}_{\psi_{\frac{2}{p}}}
    \leq p \norm{(|A_{i}x| - |A_{i}y|) \sqrt{|A_{i}x|^{2p-2}+|A_{i}y|^{2p-2}}}_{\psi_{\frac{2}{p}}}.
\end{equation*}
Also, by Hölder's inequality, we get $||XY||_{\psi_{\frac{2}{p}}} \leq ||X||_{\psi_{\frac{2r}{p}}} ||Y||_{\psi_{\frac{2s}{p}}}$ if $\frac{1}{r} + \frac{1}{s} = 1$, so it follows that
\begin{align*}
    \norm{(|A_{i}x| - |A_{i}y|) \sqrt{|A_{i}x|^{2p-2}+|A_{i}y|^{2p-2}}}_{\psi_{\frac{2}{p}}}
    &\leq \norm{|A_{i}x| - |A_{i}y|}_{\psi_{2}} \\
    &\norm{ \sqrt{|A_{i}x|^{2p-2}+|A_{i}y|^{2p-2}}}_{\psi_{\frac{2}{p-1}}}.
\end{align*}
Applying \cite[Proposition 2.6.1]{V} and Lemma \ref{cs} gives
\begin{align} \label{AA1}
    &\norm{|A_{i}x| - |A_{i}y|}_{\psi_{2}} \leq \norm{A_{i}(x-y)}_{\psi_{2}} 
    \leq
    \begin{cases}
        C_{p}(K/||A_{1,1}||_{L^{p}}) ||x-y||,&\text{ if } 1\leq p\leq 2\\
        C_{p}K||x-y||, &\text{ if } 2< p<\infty
    \end{cases}
\end{align}
and
\begin{align} \label{AA2}
    \norm{ \sqrt{|A_{i}x|^{2p-2}+|A_{i}y|^{2p-2}}}_{\psi_{\frac{2}{p-1}}}
    &\leq \bkt{C_{p}K^{2(p-1)}||x||_{2}^{2p-2}+C_{p}K^{2(p-1)}||y||_{2}^{2p-2}}^{\frac{1}{2}} \nonumber\\
    &\leq
    \begin{cases}
        C_{p}(K/||A_{1,1}||_{p})^{p-1}, &\text{ if } 1\leq p\leq 2\\
        C_{p}K^{p-1}, &\text{ if } 2<p<\infty.
    \end{cases}
\end{align}
Thus, combining (\ref{AA1}) and (\ref{AA2}) yields (\ref{Final}). Therefore, by (\ref{A1_New}) and (\ref{A2_New}), we establish (\ref{proof_case2}) when $s < 2m^{\frac{1}{p}}$, which completes the proof of Lemma \ref{v2case2}.
\end{proof}

\subsubsection{Case 3: General Vectors $x,y \in \mathbb{R}^{n}$}
\begin{lemma} \label{v2case3}
Under assumption of Theorem \ref{main}, we have
\begin{align} \label{v2case3_ineq}
    &\norm{R_{x} - R_{y}}_{\psi_{2}} \nonumber\\
    &\leq
    \begin{cases}
    C_{p} K^{p} ||A_{1,1}||_{L^{p}}^{-(p+2)} \emph{Lip}(||\cdot||_{p})||x-y||_{2}, &\text{ if } 1\leq p\leq 2\\
    C_{p} K^{p+2}   \emph{Lip}(||\cdot||_{p})||x-y||_{2}, &\text{ if } 2< p<\infty 
    \end{cases}
    \quad \forall \;x,y \in \mathbb{R}^{n}.
\end{align}
\end{lemma}
\begin{proof} Without loss of generality, we may suppose that $||x|| = 1$ and $||y|| > 1$. Set $\overline{y} \equiv \frac{y}{||y||}$.
Observe that
\begin{align}\label{RxRy}
    &\norm{R_{x} - R_{y}}_{\psi_{2}}\leq
    \norm{R_{x} - R_{\overline{y}}}_{\psi_{2}}
    +\norm{R_{\overline{y}} - R_{y}}_{\psi_{2}} \nonumber\\
    &\leq
    \begin{cases}
        \norm{R_{x} - R_{\overline{y}}}_{\psi_{2}}
        +\norm{R_{\overline{y}}}_{\psi_{2}}||y-\overline{y}||_{2}, &\text{ if } p\in [1,2]\\
        \norm{R_{x} - R_{\overline{y}}}_{\psi_{2}}
        +\norm{R_{\overline{y}}}_{\psi_{2}}C_{p}K||y-\overline{y}||_{2}, &\text{ if } p\in (2,\infty).
    \end{cases}
\end{align}
Hence, it suffices to show the reverse triangle inequality:
\begin{equation} \label{norm_triangle}
||x-\overline{y}||_{2} + ||y - \overline{y}||_{2} \leq 
\begin{cases}
    C_{p}||A_{1,1}||_{L^{p}}^{-1} ||x-y||_{2}, &\text{ if } 1\leq p\leq 2\\
    C_{p} K||x-y||_{2}, &\text{ if } 2<p<\infty
\end{cases}
\quad \forall \; ||x|| = 1, \; ||y||>1
\end{equation}
since applying Lemma \ref{v2case1}, Lemma \ref{v2case2}, and (\ref{norm_triangle}) to (\ref{RxRy}) yields Lemma \ref{v2case3}. Let $\theta$ be the angle between $x-\overline{y}$ and $y - \overline{y}$ such that $0\leq \theta \leq \pi$, i.e., $\cos \theta = \frac{\langle x-\overline{y},y-\overline{y} \rangle}{||x-\overline{y}||_{2} ||y - \overline{y}||_{2}}$. 
It is easy to see that (\ref{norm_triangle}) holds if $\frac{\pi}{2} \leq \theta \leq \pi$. Indeed, since $\cos\theta \leq 0$, applying the law of cosines gives
\begin{align*}
   (||x-\overline{y}||_{2} + ||\overline{y}-y||_{2})^{2}
   &\leq 2(||x-\overline{y}||_{2}^{2} + ||\overline{y}-y||_{2}^{2}) -4\cos(\theta) ||x-\overline{y}||_{2}||\overline{y}-y||_{2}\\
   &= 2||x-y||_{2}^{2}.
\end{align*}
In addition, if $\theta = 0$, then $\overline{y} = x$ and so there is nothing to prove.

Now, it remains to consider the case of $0 <  \theta < \frac{\pi}{2}$. Note that there are two possible positions for $y$ (as shown in Figure \ref{fig:M1}):
\begin{enumerate}
    \item 
    If $y = y_{1}$ (see the left of Figure \ref{fig:M1}), then 
    \begin{equation} \label{Case1}
    ||\overline{y} - y||_{2} + ||x-\overline{y}||_{2} \leq \frac{ \cos \theta }{\sin \theta} \sin \widetilde{\theta}||x-y||_{2} + \frac{1}{\sin \theta}\sin \widetilde{\theta} ||x-y||_{2} \leq \frac{2}{\sin \theta} ||x-y||_{2};
    \end{equation}
    \item If $y = y_{2}$ (see the right of Figure \ref{fig:M1}), then 
    \begin{align} \label{Case2}
    ||\overline{y} - y||_{2} + ||x-\overline{y}||_{2}
    &= \frac{\cos \theta}{\sin \theta}\sin \widetilde{\theta}  ||x-y||_{2} + \cos \widetilde{\theta} ||x-y||_{2} \nonumber\\
    &+ \frac{1}{\sin \theta} \sin \widetilde{\theta}||x-y||_{2} 
    \leq \frac{3}{\sin \theta} ||x-y||_{2}.
    \end{align}
\end{enumerate}
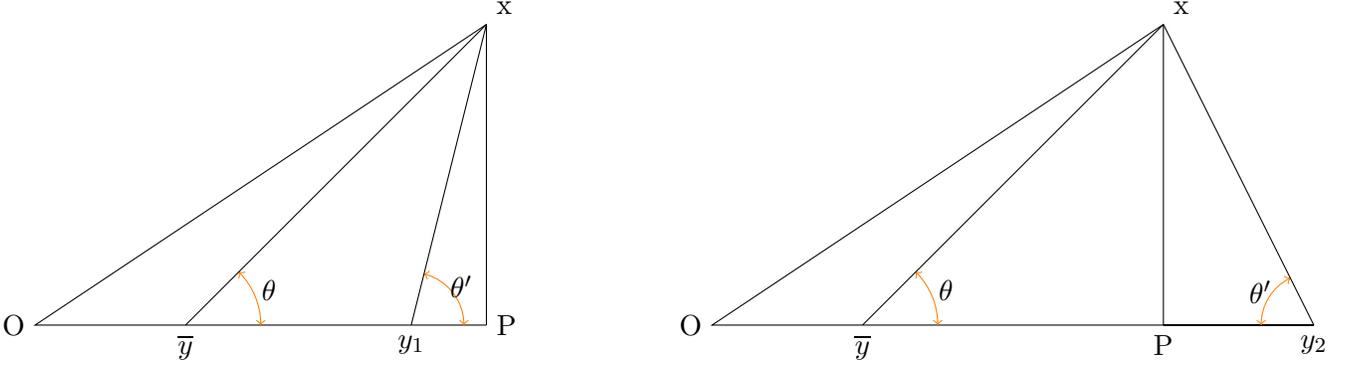
\begin{figure}
\centering
\begin{tikzpicture} [scale = 0.8]
    \draw
        (6,0) coordinate (p) node[below] {}
        -- (5,0) coordinate (y) node[below] {$y_{1}$}
        -- (2,0) coordinate (z) node[below] {$\overline{y}$}
        -- (0,0) coordinate (b) node[below] {O}
        -- (6,4) coordinate (x) node[above right] {x}
        pic["$\theta$", draw=orange, <->, angle eccentricity=1.2, angle radius=1cm]
        {angle=y--z--x}
        pic["$\widetilde{\theta}$", draw=orange, <->, angle eccentricity=1.2, angle radius=0.7cm]
        {angle=p--y--x};
    \draw
        (2,0) coordinate -- (6,4) coordinate;
    \draw
        (5,0) coordinate -- (6,4) coordinate;
    \draw
        (6,0) coordinate -- (6,4) coordinate;

    \draw
        (13,0) coordinate (p') node[below] {}
        -- (15,0) coordinate (y') node[below] {$y_{2}$}
        -- (9,0) coordinate (z') node[below] {$\overline{y}$}
        -- (7,0) coordinate (b') node[below] {O}
        -- (13,4) coordinate (x') node[above right] {x}
        pic["$\theta$", draw=orange, <->, angle eccentricity=1.2, angle radius=1cm]
        {angle=y'--z'--x'}
        pic["$\widetilde{\theta}$", draw=orange, <->, angle eccentricity=1.2, angle radius=0.7cm]
        {angle=x'--y'--p'};
    \draw
        (9,0) coordinate -- (13,4) coordinate;
    \draw
        (15,0) coordinate -- (13,4) coordinate;
    \draw
        (13,0) coordinate -- (13,4) coordinate;
\end{tikzpicture}
\caption{Left: $y = y_{1}$. Right: $y = y_{2}$.} \label{fig:M1}
\end{figure}
Thus, it suffices to show that
\begin{equation} \label{sin}
\sin \theta \geq
\begin{cases}
    C_{p} ||A_{1,1}||_{L^{p}}, &\text{ if } 1\leq p\leq 2\\
    C_{p} K^{-1}, &\text{ if } 2<p<\infty
\end{cases}
\quad \forall \; ||x|| = 1, \; ||y|| > 1 \text{ such that } 0< \theta \leq \frac{\pi}{2}.
\end{equation}
Define $B \equiv \{z \in \mathbb{R}^{n} : ||z|| \leq 1\}$ and $B_{2}(a,r) \equiv \{z \in \mathbb{R}^{n} : ||z-a||_{2} \leq r\}$. Applying (\ref{norm}) and (\ref{norm2}) yields that
$||z|| \leq \frac{1}{R_{p}} ||z||_{2}$ for every $z \in \mathbb{R}^{n}$, where $R_{p} = 1$ if $1\leq p \leq 2$; $R_{p} = \frac{1}{C_{p} K}$ if $2\leq p< \infty$. Thus, it follows that $B_{2}(0,R_{p}) \subseteq B$, $||x||_{2} > R_{p}/2$, and $||\overline{y}||_{2} > R_{p}/2$. Hence, there exists an unique $w \in \partial B_{2}(0,R_{p}/2)$ such that $\overline{Ow} \; \bot \; \overline{w \overline{y}}$. Let $\theta'$ be the angle between $\overline{x' \overline{y}}$ and $\overline{\overline{y} y}$ on the left of Figure \ref{fig:M4}. Observe that $0 < \theta' \leq \theta$. Indeed, if $\theta < \theta'$ (see the right of Figure \ref{fig:M4}), there exists $z \in B$ such that $z = r \overline{y}$ for some $r > 1$ since $B$ is convex. However, since $w \in B$, we have $||z||\leq 1$, so we get a contradiction. Therefore applying (\ref{norm}) and (\ref{norm2}) implies
\begin{equation*}
    \sin \theta \geq \sin \theta' = \frac{\overline{Ow}}{\overline{O\overline{y}}} =
    \frac{R_{p}/2}{||\overline{y}||_{2}} \geq
    \begin{cases}
        C_{p} ||A_{1,1}||_{L^{p}}, &\text{ if } 1\leq p\leq 2\\
        C_{p} K^{-1}, &\text{ if } 2<p<\infty.
    \end{cases}
\end{equation*}
Similarly, if $n\geq 3$, we consider the two dimensional space spanned by $x, \overline{y}$. Hence, (\ref{sin}) still holds, so we complete the proof of Lemma \ref{v2case3}.
\begin{figure}
\centering
\begin{tikzpicture} [scale = 0.5]
    \draw
        (1.8,-2.4) coordinate (w) node[below] {w} -- (1.8,0) coordinate;
    \draw
        (0,0) coordinate -- (1.8,-2.4) coordinate -- (9,3) coordinate (x') node[below] {x'};
    \draw
        (5,0) coordinate -- (8,4) coordinate (x) node[below] {x};
    \draw
        (0,0) coordinate (O) node[below] {O}
        -- (5,0) coordinate (z) node[below] {$\overline{y}$}
        -- (9,0) coordinate (y) node[below] {y}
        pic["$\theta$", draw=orange, <->, angle eccentricity=1.2, angle radius=1.3cm]
        {angle=y--z--x}
        pic["$\theta'$", draw=orange, <->, angle eccentricity=1.2, angle radius=1cm]
        {angle=y--z--x'};
    \draw
        (0,0) circle [radius=3];
    \draw [dashed] (5,0) -- (2,-4);

    \draw
        (14.8,-2.4) coordinate (w) node[below] {w} -- (14.8,0) coordinate;
    \draw
        (13,0) coordinate -- (14.8,-2.4) coordinate -- (22,3) coordinate (x') node[below] {x'};
    \draw
        (18,0) coordinate -- (22,1) coordinate (x'') node[below] {x};
    \draw
        (13,0) coordinate (O) node[below] {O}
        -- (18,0) coordinate (z) node[below] {$\overline{y}$}
        -- (22,0) coordinate (y) node[below] {y}
        pic["$\theta'$", draw=orange, <->, angle eccentricity=1.2, angle radius=0.7cm]
        {angle=y--z--x'}
        pic["$\theta$", draw=orange, <->, angle eccentricity=1.2, angle radius=1.2cm]
        {angle=y--z--x''};
    \draw
        (13,0) circle [radius=3];
    \draw
        [dashed] (18,0) -- (12,-2);
    \draw
        [dashed] (22,1) -- (19.88,0) coordinate (z) node[below] {z} -- (14.8,-2.4);
\end{tikzpicture}

\caption{Left: $\theta \geq \theta'$. Right: $\theta' \geq \theta$} \label{fig:M4}
\end{figure}

\end{proof}

\subsection{Proof of Corollary \ref{JL}}
Let $S = \{\frac{x-y}{||x-y||} : x,y \in T \text{ and } x\neq y\}$ and $\widetilde{S} = \{\frac{x-y}{||x-y||_{2}} : x,y \in T \text{ and } x\neq y\}$. Recall that $d_{p} ||z||_{2} \leq ||z|| \leq D_{p}||z||_{2}$ for every $z \in \mathbb{R}^{n}$, where $d_{p}$ and $D_{p}$ are defined in (\ref{d_and_D}). Thus, we have $\gamma(S) \leq \frac{1}{d_{p}} \gamma(\widetilde{S})$, so it follows that $\gamma(S) \leq \frac{1}{d_{p}} \sqrt{\log N}$ by using \cite[(9.13)]{V}. Therefore applying Theorem \ref{main} and (\ref{Lip}) gives
\begin{align*}
    &\norm{\sup_{x,y \in T, \; x\neq y}\abs{\frac{1}{m^{\frac{1}{p}}}\frac{||A(x-y)||_{p}}{||x-y||} - 1}\;}_{\psi_{2}}
    = \frac{1}{m^{\frac{1}{p}}}
    \norm{\sup_{z \in S}|R_{z}|}_{\psi_{2}} \\
    &\leq  
    \begin{cases}
    C_{p} K^{p} ||A_{1,1}||_{p}^{-(p+3)} \sqrt{\log(N)}  \frac{1}{m^{1/2}},& \text{ if } 1\leq p\leq 2\\
    C_{p} K^{p+2} \sqrt{\log(N)} \frac{1}{m^{1/p}},& \text{ if } 2< p<\infty,
    \end{cases}
\end{align*}
which implies (\ref{JL_EQ}).

\bibliographystyle{plainurl}
\bibliography{mybibliography.bib}

\begin{thebibliography}{10}

\bibitem{intro1}
Bo~Brinkman and Moses Charikar.
\newblock On the impossibility of dimension reduction in $\ell_{1}$.
\newblock {\em J. ACM}, 52(5):766–788, 2005.
\newblock \href {https://doi.org/10.1145/1089023.1089026}
  {\path{doi:10.1145/1089023.1089026}}.

\bibitem{Chen}
Wei-Kuo Chen and Arnab Sen.
\newblock {On $\ell_{p}$-Gaussian–Grothendieck Problem}.
\newblock {\em International Mathematics Research Notices}, 11 2021.
\newblock \href {https://doi.org/10.1093/imrn/rnab311}
  {\path{doi:10.1093/imrn/rnab311}}.

\bibitem{chow1997probability}
Y.S. Chow, H.~Teicher, G.~Casella, S.~Fienberg, and I.~Olkin.
\newblock {\em Probability Theory: Independence, Interchangeability,
  Martingales}.
\newblock Springer Texts in Statistics. Springer New York, 1997.
\newblock \href {https://doi.org/10.1007/978-1-4612-1950-7}
  {\path{doi:10.1007/978-1-4612-1950-7}}.

\bibitem{a_sub_exponential}
Friedrich Götze, Holger Sambale, and Arthur Sinulis.
\newblock {Concentration inequalities for polynomials in
  $\alpha$-sub-exponential random variables}.
\newblock {\em Electronic Journal of Probability}, 26:1 -- 22, 2021.
\newblock \href {https://doi.org/10.1214/21-EJP606}
  {\path{doi:10.1214/21-EJP606}}.

\bibitem{diamond}
{James R.} Lee and Assaf Naor.
\newblock Embedding the diamond graph in $\emph{L}_{p}$ and dimension reduction
  in $\emph{L}_{1}$.
\newblock {\em Geometric and Functional Analysis}, 14(4):745--747, 2004.
\newblock \href {https://doi.org/10.1007/s00039-004-0473-8}
  {\path{doi:10.1007/s00039-004-0473-8}}.

\bibitem{intro2}
Ping Li.
\newblock Stable random projections and conditional random sampling, two
  sampling techniques for modern massive datasets.
\newblock 2007.

\bibitem{V2}
Christopher Liaw, Abbas Mehrabian, Yaniv Plan, and Roman Vershynin.
\newblock A simple tool for bounding the deviation of random matrices on
  geometric sets.
\newblock In {\em Geometric Aspects of Functional Analysis: Israel Seminar},
  pages 277--299. Springer International Publishing, 2017.
\newblock \href {https://doi.org/10.1007/978-3-319-45282-1_18}
  {\path{doi:10.1007/978-3-319-45282-1_18}}.

\bibitem{Tstable}
Michael~B. Marcus and Michel Talagrand.
\newblock Chevet{\textquotesingle}s theorem for stable processes {II}.
\newblock {\em Journal of Theoretical Probability}, 1(1):65--92, January 1988.
\newblock \href {https://doi.org/10.1007/bf01076288}
  {\path{doi:10.1007/bf01076288}}.

\bibitem{sJ3}
Shahar Mendelson.
\newblock On weakly bounded empirical processes.
\newblock {\em Mathematische Annalen}, 340(1):293--314, 2008.
\newblock \href {https://doi.org/10.1007/s00208-007-0152-9}
  {\path{doi:10.1007/s00208-007-0152-9}}.

\bibitem{sJ2}
Shahar Mendelson, Alain Pajor, and Nicole Tomczak-Jaegermann.
\newblock Uniform uncertainty principle for bernoulli and subgaussian
  ensembles.
\newblock {\em Constructive Approximation}, 28:277--289, 12 2008.
\newblock \href {https://doi.org/10.1007/s00365-007-9005-8}
  {\path{doi:10.1007/s00365-007-9005-8}}.

\bibitem{sJ1}
Shahar Mendelson and Nicole Tomczak-Jaegermann.
\newblock A subgaussian embedding theorem.
\newblock {\em Israel Journal of Mathematics}, 164:349--364, 2008.
\newblock \href {https://doi.org/10.1007/s11856-008-0034-1}
  {\path{doi:10.1007/s11856-008-0034-1}}.

\bibitem{alpha'}
Holger Sambale.
\newblock Some notes on concentration for $\alpha$-subexponential random
  variables.
\newblock 2020.
\newblock \href {https://doi.org/10.48550/ARXIV.2002.10761}
  {\path{doi:10.48550/ARXIV.2002.10761}}.

\bibitem{Schechtman}
Gideon Schechtman.
\newblock Two observations regarding embedding subsets of euclidean spaces in
  normed spaces.
\newblock {\em Advances in Mathematics}, 200(1):125--135, 2006.
\newblock \href {https://doi.org/10.1016/j.aim.2004.11.003}
  {\path{doi:10.1016/j.aim.2004.11.003}}.

\bibitem{stroock}
D.W. Stroock.
\newblock {\em Essentials of Integration Theory for Analysis}.
\newblock Graduate Texts in Mathematics. Springer New York, 2011.
\newblock \href {https://doi.org/10.1007/978-1-4614-1135-2}
  {\path{doi:10.1007/978-1-4614-1135-2}}.

\bibitem{V}
Roman Vershynin.
\newblock {\em High-Dimensional Probability: An Introduction with Applications
  in Data Science}.
\newblock Cambridge Series in Statistical and Probabilistic Mathematics.
  Cambridge University Press, 2018.
\newblock \href {https://doi.org/10.1017/9781108231596}
  {\path{doi:10.1017/9781108231596}}.

\end{thebibliography}

\end{document}